\documentclass[12pt]{article}
\usepackage[cp1251]{inputenc}
\usepackage[russian]{babel}
\usepackage[T2B]{fontenc}
\usepackage{indentfirst}
\usepackage[pdftex, colorlinks, urlcolor=blue, pdfborder={0 0 0 [0 0]}]{hyperref}
\usepackage{amsfonts,amsmath,amssymb,theorem}
\usepackage{array}

\newtheorem{Theorem}{Теорема}

\newtheorem{rem}{Замечание}
\newtheorem{dfn}{Определение}
\newtheorem{exm}{Пример}
\newenvironment{proof}
{\vspace{1pt}\par{\sl%
Доказательство.\,\ }}%
{\noindent\vspace{1pt}}

\begin{document}
\title
{Построение программных управлений с вероятностью 1 для динамической системы с пуассоновскими возмущениями}
\author{Е. В. Карачанская \\
{Тихоокеанский государственный университет, Россия, Хабаровск}}

\date{}
\maketitle

\begin{abstract}
В статье предлагается метод построения программных управлений с вероятностью 1 для динамической системы, подверженной пуассоновским возмущениям.
\end{abstract}

\section*{Введение}
Одной из задач управления является организация управления
динамической системой таким образом, чтобы при ее эволюции важные
характеристики системы (в том числе и зависящие от положения
системы), сохранялись. В реальном пространстве на динамическую
систему оказывают влияние случайные факторы. Наиболее удачно это
случайное воздействие можно описать с помощью винеровских и
пуассоновских процессов. В последнее время было предложено
несколько способов построения управлений системой при наличии
винеровских процессов. В данной работе предлагается метод
построения программного управления с вероятностью 1 для системы,
подверженной возмущениям в виде винеровских и пуассоновских
процессов. Предложенный метод основан на понятии первого интеграла
для системы стохастических дифференциальных уравнений с
винеровскими и пуассоновскими возмущениями
{\rm{\cite{D-78,D-89,D-02}}} и алгоритме построения автоморфной
функции {\rm{\cite{D-03}}}.

\section{Стохастический первый интеграл системы СДУ}

Пусть $\bf{x}\in \mathbb{R}^{n}$, $\bf{x}(t)$ -- случайных
процесс, являющийся решением системы стохастических
дифференциальных уравнений
\begin{equation}\label{1}
\begin{array}{l}
  dx_{i}(t)=a_{i}(t;{\bf{x}}(t))dt+\displaystyle\sum\limits_{k=1}^{m}b_{ik}(t;{\bf{x}}(t))
  dw_{k}(t)+
  \displaystyle\int\limits_{R(\gamma)}g_{i}(t;{\bf{x}}(t);\gamma)\nu(dt;d\gamma)\\
{\bf{x}}(t)={\bf{x}}(t,{\bf{x}}_{o})\Bigl|_{t=0}={\bf{x}}_{0},\ \
\ i=\overline{1,n}, \ \ t\geq 0,
\end{array}
\end{equation}
где ${\bf{w}}(t))$ -- $m$-мерный винеровский процесс,
$\nu(t;\Delta\gamma)$ -- однородная по $t$, нецентрированная мера
Пуассона {\rm{\cite{GS-68}}}. Эту систему можно
записать в векторной форме
\begin{equation*}
d {\bf x}(t)= A(t;{\bf x}(t)) dt +  B(t;{\bf x}(t))  d {\bf
w}(t)+\displaystyle\int\limits_{R(\gamma)}\nu(dt;d\gamma)\cdot
G(t;{\bf x}(t);\gamma).
\end{equation*}
Относительно коэффициентов $a_{i}(t;{\bf{x}})$,
$b_{ik}(t;{\bf{x}})$ и $g_{i}(t;{\bf{x}};\gamma)$ уравнения
(\ref{1}) будем предполагать, что они выбраны таким образом, чтобы
были обеспечены условия существования и единственности решения,
как и во всех уравнениях, рассматриваемых ниже.

В \cite{D-78} было введено понятие первого интеграла для системы
стохастических дифференциальных уравнений Ито (без пуассоновской
составляющей), в \cite{D-02} -- понятие стохастического
первого интеграла для системы обобщенных стохастических
дифференциальных уравнений Ито с центрированной пуассоновской
мерой. Введем аналогичное понятие для случая наличия
нецентрированной меры Пуассона.

Пусть $u(t;{\bf x};\omega)$ -- случайная функция, определенная на
том же вероятностном пространстве, что и решение системы
(\ref{1}).
\begin{dfn}\label{df1}
Случайную функцию $u(t;{\bf x};\omega)$ назовем стохастическим
первым интегралом системы  (\ref{1}), если с вероятностью 1
$$
u(t;{\bf x}(t; \mathbf{x}(0));\omega)=u(0;\mathbf{x}(0))
$$
для любого решения $\mathbf{x}(t;\mathbf{x}(0);\omega)$ системы
(\ref{1}).
\end{dfn}

Чтобы функция $u(t;{\bf x};\omega)$ была первым интегралом системы
(\ref{1}), должны выполняться условия $\left.\mathcal{L}\right)$:
\begin{enumerate}
    \item $\displaystyle\sum\limits_{i=1}^{n}b_{i\,k}(t;{\bf x})\frac{\partial u(t;{\bf x})}{\partial
x_{i}}=0$, для всех $k=\overline{1,m}$ (компенсация винеровского
возмущения);
    \item $\displaystyle\frac{\partial u(t;{\bf x})}{\partial
t}+\displaystyle\sum\limits_{i=1}^{n}\frac{\partial u(t;{\bf
x})}{\partial x_{i}}\Bigl[a_{i}(t;{\bf x})-
\displaystyle\frac{1}{2}\sum\limits_{k=1}^{m}\sum\limits_{j=1}^{n}\,b_{j\,k}(t;{\bf
x})\frac{\partial b_{i\,k}(t;{\bf x})}{\partial x_{j}}\Bigr]=0$
(независимость от времени);
    \item $u(t;{\bf x})-u\Bigl(t;{\bf x}+G(t;{\bf
    x};\gamma)\Bigr)=0$ для любых $\gamma\in R(\gamma)$ во всей области
  определения процесса (компенсация пуассоновских скачков).
\end{enumerate}
\begin{rem}
В случае, когда рассматриваем конкретную реализацию, т. е.
параметр $\omega$ в дальнейшем не влияет, неслучайную функцию
$u(t;{\bf x})$ можно считать детерминированным первым интегралом
стохастической системы.
\end{rem}

В \cite{D-02} было введено понятие стохастического первого
интеграла для центрированной пуассоновской меры, и полученные
условия для его существования учитывают необходимость задания
плотности интенсивности пуассоновского распределения в отличие от
предложенного в данной статье. Таким образом, безразлично, каков
вероятностный закон имеют интенсивности пуассоновских скачков.

\section{Построение системы  обобщенных СДУ с заданным первым интегралом}
Определим вид системы обобщенных стохастических уравнений Ито с
начальными данными, имеющей известный стохастический первый
интеграл.

\begin{Theorem}\label{teor-1}
Пусть функция $u(t,{\bf x})$ -- непрерывна вместе со своими
производными по совокупности переменных $(t,{\bf x})$ и случайная
функция $u(t,{\bf x};\omega)$ определена на том же вероятностном
пространстве, что и решение системы стохастических
дифференциальных уравнений
\begin{equation}\label{1m}
\begin{array}{l}
  d{\bf{x}}(t)=A(t;{\bf{x}}(t))dt+B(t;{\bf{x}}(t))
  d\mathbf{w}(t)+
  \displaystyle\int\limits_{R(\gamma)}\nu(dt;d\gamma)\cdot
G(t;{\bf x}(t);\gamma)\\
{\bf{x}}(t)={\bf{x}}(t,{\bf{x}}_{o})\Bigl|_{t=0}={\bf{x}}_{0},\ \
\ t\geq 0,
\end{array}
\end{equation}
где ${\bf{x}}\in \mathbb{R}^{n}$, $n\geq 2$;  ${\bf{w}}(t)$ --
$m$-мерный винеровский процесс; $\nu(t;\Delta\gamma)$ --
однородная по $t$, нецентрированная мера Пуассона. Если $u(t,{\bf
x};\omega)$ является стохастическим первым интегралом системы
(\ref{1m}), то коэффициенты уравнения (\ref{1m}) и функция
$u(t,{\bf x})$ связаны следующими соотношениями:
\begin{enumerate}
    \item коэффициенты
$B_{k}(t;{\bf{x}})=\displaystyle\sum\limits_{i=1}^{n}b_{ik}(t;{\bf{x}})\vec{e}_{i}$
    $(k=\overline{1,m})$
    -- столбцы матрицы $B(t;{\bf{x}})$, $B_{k}(t;{\bf{x}})\in \{q_{00}\cdot M_{n+1,0}\}$, $M_{n+1,0}$ -- минор элемента $h_{n+1,0}$ матрицы $H(t;{\bf{x}})$:
  \begin{equation}\label{B}
H(t;{\bf{x}})=\left(
\begin{array}{cccc}
  \vec{e}_{0} & \vec{e}_{1} &\ldots & \vec{e}_{n} \\
  \displaystyle\frac{\partial u(t; {\bf{x}})}{\partial t}  & \displaystyle\frac{\partial u(t; {\bf{x}})}{\partial x_{1}}  & \ldots &
  \displaystyle\frac{\partial u(t; {\bf{x}})}{\partial x_{n}}  \\
  h_{30} & h_{31} &\ldots & h_{3n}\\
  \ldots & \ldots &\ldots & \ldots  \\
  h_{n+1,0} &  h_{n+1,1} & \ldots & h_{n+1,n}
\end{array}
\right),
\end{equation}
    \item коэффициент
    $A(t;{\bf{x}})$
    принадлежит множеству функций, определяемых условием
\begin{equation}\label{A}
{{A}}(t;{\bf x})\in \left\{ R(t;{\bf x}) + \displaystyle
\frac{1}{2}\, \sum_{k=1}^{n}\displaystyle \biggl[\frac{\partial
B_{k}(t;{\bf x})}{\partial {\bf x}}\biggr]\cdot B_{k}(t;{\bf
x})\right\},
\end{equation}
где
$\displaystyle \biggl[\frac{\partial B_{k}(t;{\bf x})}{\partial
{\bf x}}\biggr]$ -- матрица Якоби для векторной функции
$B_{k}(t;{\bf x})$;
$C(t;{\bf x})$ -- алгебраическое дополнение элемента $\vec{e}_{0}$ матрицы $H(t;{\bf{x}})$ и $\det C(t;{\bf x})\neq 0$; матрица-столбец $ R(t;{\bf x})$, компоненты которой $r_{i}(t;{\bf
x}) $, $i=\overline{1,n}$, определяются следующим образом:
$$ C^{-1}(t;{\bf x}) \cdot \det H(t;{\bf{x}})= \vec{e}_{0}+
\displaystyle \sum\limits_{i=1}^{n} r_{i}(t;{\bf x})\vec{e}_{i};
$$
    \item коэффициент $G(t;{\bf x};\gamma)=
    \displaystyle\sum\limits_{i=1}^{n}g_{i}(t;{\bf{x}};\gamma)\vec{e}_{i}$
     при пуассоновской мере определяется представлением
     $G(t;{\bf x};\gamma)=\mathbf{y}(t;{\bf x};\gamma)-{\bf x}$,
где $\mathbf{y}(t;{\bf x};\gamma)$ -- решение системы
дифференциальных уравнений
\begin{equation}\label{Y}
\displaystyle\frac{\partial \mathbf{y}(\cdot;\gamma)}{\partial
\gamma}= \det\left[
\begin{array}{cccc}
  \vec{e}_{1}& \vec{e}_{2} & \cdots & \vec{e}_{n} \\
  \displaystyle\frac{\partial u(t; \mathbf{y}(\cdot;\gamma))}{\partial y_{1}} &
  \displaystyle\frac{\partial u(t; \mathbf{y}(\cdot;\gamma))}{\partial y_{2}} & \cdots &
  \displaystyle\frac{\partial u(t; \mathbf{y}(\cdot;\gamma))}{\partial y_{n}} \\
  \varphi_{31}(t;\mathbf{y}(\cdot;\gamma)) &
  \varphi_{32}(t;\mathbf{y}(\cdot;\gamma)) & \cdots &
  \varphi_{3n}(t;\mathbf{y}(\cdot;\gamma)) \\
  \cdots & \cdots & \cdots & \cdots \\
   \varphi_{n1}(t;\mathbf{y}(\cdot;\gamma)) &
  \varphi_{n2}(t;\mathbf{y}(\cdot;\gamma)) & \cdots &
  \varphi_{nn}(t;\mathbf{y}(\cdot;\gamma))
\end{array}
\right]
\end{equation}
удовлетворяющее начальному условию
$\mathbf{y}(t;\mathbf{x};\gamma)\Bigl|_{\gamma=0}=\mathbf{x}$.
\end{enumerate}
Относительно произвольных функций 
$h_{ij}=h_{ij}(t,\mathbf{x} )$,
$\varphi_{ij}=\varphi_{ij}(t;\mathbf{y}(\cdot;\gamma))$ полагаем,
что они выбраны таким образом, чтобы каждое семейство функций
 $\Bigl\{h_{i}\Bigr\}$,
$\Bigl\{\varphi_{i}\Bigr\}$, определяемое условиями:
$$
h_{ij}(t,\mathbf{x} )=\displaystyle\frac{\partial
h_{i}(t,\mathbf{x} )}{\partial x_{j}}, \ \ \ \ \
\varphi_{ij}(t;\mathbf{y}(\cdot;\gamma))=\displaystyle\frac{\partial
\varphi_{i}(t;\mathbf{y}(\cdot;\gamma))}{\partial y_{j}},
$$
составляло вместе с функцией $u(t;\mathbf{x})$ совокупность
независимых функций.
\end{Theorem}
\begin{proof}
Доказательство состоит из 3-х частей.

${\bf 1}.$ Воспользуемся первым из условий
$\left.\mathcal{L}\right)$:
$\displaystyle\sum\limits_{i=1}^{n}b_{i\,k}(t;{\bf
x})\displaystyle\frac{\partial u(t;{\bf x})}{\partial x_{i}}=0$,
для всех $k=\overline{1,m}$. Если
$B_{k}(t;{\bf{x}})=\displaystyle\sum\limits_{i=1}^{n}b_{ik}(t;{\bf{x}})\vec{e}_{i}$
и $\nabla_{\mathbf{x}}u(t;{\bf
x})=\displaystyle\sum\limits_{i=1}^{n}\frac{\partial u(t;{\bf
x})}{\partial x_{i}}\, \vec{e}_{i} $, то это условие -- есть
условие ортогональности векторов
 $B_{k}(t;{\bf{x}})$ и $\nabla_{\mathbf{x}}u(t;{\bf
x})$. Опираясь на определение векторного произведения в
пространстве $\mathbb{R}^{n}$ и его свойства, получаем утверждение
для коэффициентов $B_{k}(t;{\bf{x}})$ и, соответственно, матрицы
$B(\cdot)=\Bigl(B_{1}(\cdot),\ldots,B_{m}(\cdot)\Bigr)$:
$$
B_{k}(t;{\bf{x}})\in \left\{
\begin{array}{c}
q_{00}\cdot \det \left(
\begin{array}{ccc}
  \vec{e}_{1} & \ldots & \vec{e}_{n} \\
  \displaystyle\frac{\partial u(t; {\bf{x}})}{\partial x_{1}}  & \ldots &
  \displaystyle\frac{\partial u(t; {\bf{x}})}{\partial x_{n}}  \\
  f_{31} & \ldots & f_{3n}\\
  \ldots & \ldots & \ldots  \\
  f_{n1} & \ldots & f_{nn}
\end{array}
\right)
\end{array}
\right\};
$$
где функции $f_{i}=f_{i}(t,\mathbf{x} )$, \, $i=\overline{3,n}$,\,
такие, что $f_{ij}(t,\mathbf{x} )=\displaystyle\frac{\partial
f_{i}(t,\mathbf{x} )}{\partial x_{j}}$ вместе с
 функцией $u(t,\mathbf{x} )$ образуют совокупность независимых
 функций.

${\bf 2}.$ Воспользуемся вторым из условий
$\left.\mathcal{L}\right)$:
$$
\displaystyle\frac{\partial u(t;{\bf x})}{\partial
t}+\displaystyle\sum\limits_{i=1}^{n}\frac{\partial u(t;{\bf
x})}{\partial x_{i}}\Bigl[a_{i}(t;{\bf x})-
\displaystyle\frac{1}{2}\sum\limits_{k=1}^{m}\sum\limits_{j=1}^{n}\,b_{j\,k}(t;{\bf
x})\frac{\partial b_{i\,k}(t;{\bf x})}{\partial x_{j}}\Bigr]=0.
$$
Пусть
$$Q(t;{\bf x})=1+
\displaystyle\sum\limits_{i=1}^{n}a_{i}(t;{\bf{x}})-
\displaystyle\frac{1}{2}\sum\limits_{i=1}^{n}\sum\limits_{k=1}^{m}\sum\limits_{j=1}^{n}\,b_{j\,k}(t;{\bf
x})\frac{\partial b_{i\,k}(t;{\bf x})}{\partial x_{j}}.
$$
Следуя схеме, изложенной в работе {\rm{\cite{D-89}}}, введем в
рассмотрение векторы: обобщенный градиент
$$\square
u(t;{\bf x})= \displaystyle\frac{\partial u(t;{\bf x})}{\partial
t}\, \vec{e}_{0}+\sum\limits_{i=1}^{n}\frac{\partial u(t;{\bf
x})}{\partial x_{i}}\, \vec{e}_{i}
$$
и
$$\overrightarrow{Q}(t;{\bf x})=\vec{e}_{0}+
\displaystyle\sum\limits_{i=1}^{n}a_{i}(t;{\bf{x}})\vec{e}_{i}-
\displaystyle\frac{1}{2}\sum\limits_{i=1}^{n}\sum\limits_{k=1}^{m}\sum\limits_{j=1}^{n}\,b_{j\,k}(t;{\bf
x})\frac{\partial b_{i\,k}(t;{\bf x})}{\partial
x_{j}}\,\vec{e}_{i}.
$$
Указанное условие означает, что векторы $\square u(t;{\bf x})$ и
$\overrightarrow{Q}(t;{\bf x})$ ортогональны. Воспользовавшись
снова определением векторного произведения и его свойства,
получаем формулу \eqref{B}:
$
\overrightarrow{Q}(t;{\bf x})\in \left\{\det H \right\},
$
где функции $h_{i}=f_{i}(t,\mathbf{x} )$, \,
$i=\overline{3,n+1}$,\, такие, что $h_{ij}(t,\mathbf{x}
)=\displaystyle\frac{\partial h_{i}(t,\mathbf{x} )}{\partial
x_{j}}$ вместе с
 функцией $u(t,\mathbf{x} )$ образуют совокупность независимых
 функций. Без ограничения общности, будем считать, что $h_{ij}(t,\mathbf{x})=f_{ij}(t,\mathbf{x}.
)$
Введем вектор
$$
\overrightarrow{{\widetilde{A}}}(t;{\bf x})= \vec{e}_{o}+
\displaystyle\sum\limits_{i=1}^{n}a_{i}(t;{\bf x})\vec{e}_{i}=
\overrightarrow{Q}(t;{\bf x}) +
\displaystyle\frac{1}{2}\displaystyle\sum\limits_{i=1}^{n}\sum\limits_{k=1}^{m}\sum\limits_{j=1}^{n}\,b_{j\,k}(t;{\bf
x})\frac{\partial b_{i\,k}(t;{\bf x})}{\partial
x_{j}}\,\vec{e}_{i}
$$
 и, в силу того, что коэффициент
при $\vec{e}_{o}$ должен быть равен 1, получаем:
\begin{equation*}
\overrightarrow{{\widetilde{A}}}(t;{\bf x})\in \left\{
C^{-1}(t;{\bf x}) \cdot \det H(t;{\bf x}) + \displaystyle \frac{1}{2}\,
\displaystyle\sum\limits_{i=1}^{n}\sum\limits_{k=1}^{m}\sum\limits_{j=1}^{n}\,b_{j\,k}(t;{\bf
x})\frac{\partial b_{i\,k}(t;{\bf x})}{\partial
x_{j}}\,\vec{e}_{i}\right\},
\end{equation*}
где $C(t;{\bf x}) $ -- алгебраическое дополнение элемента
$\vec{e}_{0}$ матрицы $H(t;{\bf x})$, $\det C(t;{\bf x})\neq 0 $. Поскольку
вектор $C^{-1}(t;{\bf x}) \cdot \det H(t;{\bf x}) $ можно записать в виде:
$$ C^{-1}(t;{\bf x}) \cdot \det H= \vec{e}_{0}+
\displaystyle \sum\limits_{i=1}^{n} r_{i}(t;{\bf x})\vec{e}_{i},
$$
то введем матрицу-столбец $ R(t;{\bf x})$ с компонентами
$r_{i}(t;{\bf x}) $, $i=\overline{1,n}$.

Определение произведения матриц в данном случае допускает
представление:
\begin{equation*}
\displaystyle\sum\limits_{i=1}^{n}\sum\limits_{k=1}^{m}\sum\limits_{j=1}^{n}\,b_{j\,k}(t;{\bf
x})\frac{\partial b_{i\,k}(t;{\bf x})}{\partial
x_{j}}=\sum_{k=1}^{n}\displaystyle \biggl[\frac{\partial
B_{k}(t;{\bf x})}{\partial {\bf x}}\biggr]\cdot B_{k}(t;{\bf x}),
\end{equation*}
где $\displaystyle \biggl[\frac{\partial B_{k}(t;{\bf
x})}{\partial {\bf x}}\biggr]$ -- матрица Якоби для векторной
функции $B_{k}(t;{\bf x})$. Следовательно, 
${{A}}(t;{\bf x})$ определяется суммой матриц \eqref{A}:
\begin{equation*}
{{A}}(t;{\bf x})\in \left\{ R(t;{\bf x})+ \displaystyle
\frac{1}{2}\, \sum\limits_{k=1}^{m}\sum_{k=1}^{n}\displaystyle
\biggl[\frac{\partial B_{k}(t;{\bf x})}{\partial {\bf
x}}\biggr]\cdot B_{k}(t;{\bf x})\right\},
\end{equation*}

${\bf 3}.$ Исходя из третьего условия в
$\left.\mathcal{L}\right)$, для  любых $\gamma\in R(\gamma)$
должно выполняться условие:
$$
u(t;{\bf x};\omega)-u\Bigl(t;{\bf x}+G(t;{\bf
x};\gamma);\omega\Bigr)=0.
$$
Это означает, что функция $u(t;{\bf x};\omega)$ является
автоморфной при преобразовании ее аргумента ${\bf x}$ с помощью
функции $G(t;{\bf x};\gamma)$. Определим условия, налагаемые на
эту функцию. Следуя {\rm{\cite{D-02}}}, положим: ${\bf y}(t;{\bf
x};\gamma)={\bf x}+G(t;{\bf x};\gamma)$. Для упрощения записи
будем опускать параметр $\omega$. Тогда
$
u(t;{\bf x})=u\Bigl(t;{\bf y}(t;{\bf x};\gamma)\Bigr)
$
любых $\gamma\in R(\gamma)$ или $\displaystyle
\frac{\partial u(t;{\bf x})}{\partial \gamma}=0$ и 
$$
\displaystyle \frac{\partial u\Bigl(t;{\bf y}(t;{\bf
x};\gamma)\Bigr)}{\partial
\gamma}\equiv\sum\limits_{i=1}^{n}\frac{\partial u\Bigl(t;{\bf
y}(t;{\bf x};\gamma)\Bigr)}{\partial y_{i}}\frac{\partial
y_{i}(t;{\bf x};\gamma)}{\partial \gamma}=0.
$$
Последнее равенство означает, что векторы $\nabla_{{\bf
y}}u\Bigl(t;{\bf
y}(\cdot;\gamma)\Bigr)=\displaystyle\sum\limits_{i=1}^{n}\frac{\partial
u\Bigl(t;{\bf y}(\cdot;\gamma)\Bigr)}{\partial y_{i}}\,
\vec{e}_{i}$ и $\displaystyle\frac{\partial {\bf y}(\cdot;\gamma)
}{\partial \gamma}=\sum\limits_{i=1}^{n}\frac{\partial
y_{i}(\cdot;\gamma)}{\partial \gamma}\, \vec{e}_{i}$ ортогональны,
и следовательно, связаны соотношением:
\begin{equation}\label{Gam}
\displaystyle\frac{\partial {\bf y}(\cdot;\gamma) }{\partial
\gamma}\in \left\{\det\left[
\begin{array}{ccc}
  \vec{e}_{1} & \ldots & \vec{e}_{n} \\
  \displaystyle \frac{\partial
u\Bigl(t;{\bf y}(\cdot;\gamma)\Bigr)}{\partial y_{1}}& \ldots &
\displaystyle \frac{\partial
u\Bigl(t;{\bf y}(\cdot;\gamma)\Bigr)}{\partial y_{n}} \\
  \varphi_{31} & \ldots& \varphi_{3n} \\
\ldots & \ldots& \ldots \\
\varphi_{n1} & \ldots& \varphi_{nn}
\end{array}
\right]\right\},
\end{equation}
где функции $\varphi_{i}(t;{\bf y})$, $i=\overline{3,n}$, такие
что $\varphi_{ij} =\displaystyle \frac{\partial \varphi_{i}(t;{\bf
y})}{\partial y_{j}}$, составляют с функцией $u\Bigl(t;{\bf
y}(\cdot;\gamma)\Bigr)$ систему независимых функций. Поскольку
${\bf y}(t;{\bf x};\gamma)={\bf x}+G(t;{\bf x};\gamma)$, то
(\ref{Gam}) можно рассматривать как систему дифференциальных
уравнений, в которой неизвестной является функция ${\bf
y}(\cdot;\gamma)$. Разложим определитель (\ref{Gam}) по первой
строке. Следовательно, $\displaystyle\frac{\partial {\bf
y}(\cdot;\gamma) }{\partial \gamma}= \alpha
\displaystyle\sum\limits_{i=1}^{n}S_{i}({\bf y}(\cdot;\gamma))\,
\vec{e}_{i}$, где $\alpha$ -- произвольная функция, не зависящая
от ${\bf y}$. Таким образом, получаем систему дифференциальных
уравнений
\begin{equation*}
\left\{\begin{array}{c}
 \displaystyle\frac{\partial  y_{1}(\cdot;\gamma) }{\partial
\gamma}= \alpha S_{1}({\bf y}(\cdot;\gamma)),\\
\cdots \\
 \displaystyle\frac{\partial  y_{n}(\cdot;\gamma)
}{\partial \gamma}= \alpha S_{n}({\bf y}(\cdot;\gamma)).
\end{array}
\right.
\end{equation*}
Пусть ${\bf y}(t;{\bf x};\gamma;\theta)$ -- решение этой системы, где $\theta$ -- вектор постоянных, появившихся при
ее интегрировании. Поскольку условие (\ref{u}) должно выполняться
для любых $t$,\, ${\bf x}$ и $\gamma$, то
\begin{equation*}
u(t;{\bf x})=u\Bigl(t;{\bf y}(t;{\bf
x};\gamma_{1};\theta)\Bigr)=u\Bigl(t;{\bf x}+ G(t;{\bf
x};\gamma_{1};\theta)\Bigr)=u\Bigl(t;{\bf x}+ G(t;{\bf
x};\gamma_{2};\theta)\Bigr),
\end{equation*}
В частном случае, при некотором значении $\gamma=\gamma_{o}$
последнее равенство будет определяться условием $G(t;{\bf
x};\gamma_{o};\theta)=0$. Без ограничения общности (при отсутствии
скачка), положим:
$$
G(t;{\bf x};\gamma_{o};\theta)\equiv G(t;{\bf x};0)=0.
$$
Следовательно, функция $G(t;{\bf x};\gamma)$, обеспечивающая
автоморфизм функции $u(t;{\bf x})$ определяется представлением
     $G(t;{\bf x};\gamma)=\mathbf{y}(t;{\bf x};\gamma)-{\bf x}$,
где $\mathbf{y}(t;{\bf x};\gamma)$ -- решение системы
дифференциальных уравнений (\ref{Y}) при начальном условии
$\mathbf{y}(t;\mathbf{x};\gamma)\Bigl|_{\gamma=0}=\mathbf{x}$. Таким образом, последнее утверждение теоремы доказано.
\end{proof}

\section{Построение программных управлений}
Очень часто возникает задачи управления динамической системой, в
которой сохраняются заданные функционалы
{\rm{\cite{D_95}}}, причем влияние случайных возмущений,
действующих на данную систему, должно быть сведено к минимуму.
Понятие стохастического первого интеграла системы стохастических
дифференциальных уравнений с винеровскими и пуассоновскими
возмущениями позволяет строить такие управления с вероятностью 1,
то есть полностью исключая влияние данных случайных возмущений.

По аналогии с {\rm{\cite{Ch-U1}}} введем следующее определение.
\begin{dfn}\label{df2}
Программным движением стохастической системы
\begin{equation}\label{2}
d {\bf x}(t)= \Bigl[ P(t;{\bf x}(t)) +  Q(t;{\bf x}(t)) \cdot {\bf
s}(t;{\bf x}(t))  \Bigr] dt +  B(t;{\bf x}(t))  d {\bf
w}(t)+\displaystyle\int\limits_{R(\gamma)}G(t;{\bf
x}(t);\gamma)\nu(dt;d\gamma),
\end{equation}
где ${\bf w}(t)$ -- $m$-мерный винеровский процесс;
$\nu(t;\triangle \gamma)$ -- нецентрированная пуассоновская мера,
будем называть решение ${\bf x}(t; {\bf x}_{o},{\bf s};\omega)$,
позволяющее с вероятностью {\rm 1} при некотором управлении
(программном управлении) ${\bf s}(t;{\bf x})$ \underline{для всех
$t $ } оставаться на неслучайном интегральном многообразии $
u\Bigl(t;{\bf x}(t;{\bf x}_{o})\Bigr)= u(0;{\bf x}_{o}), $
являющимся первым интегралом уравнения (\ref{2}) при заданных
начальных условиях
$${\bf x}(t;{\bf x}_{o})\Bigr|_{t=0}={\bf x}_{o}.
$$
\end{dfn}

Таким образом, можно построить программное управление с
вероятностью 1 для динамической системы, подверженной случайному
воздействию винеровских и пуассоновских процессов.

\begin{Theorem}\label{teor-2}
Программное управление, позволяющее с вероятностью 1 динамической
системе (\ref{2}) при наличии винеровских и пуассоновских
возмущений оставаться на интегральном многообразии $ u\Bigl(t;{\bf
x}(t;{\bf x}_{o});\omega\Bigr)= u(0;{\bf x}_{o})$, является
решением системы, состоящей их уравнений (\ref{2}) и (\ref{1m}), в
которой коэффициенты второго уравнения и соответствующие
коэффициенты первого определяются в соответствии с
Теоремой~\ref{teor-1}. При этом определяются также реакции на
случайные возмущения, обеспечивающие это программное управление.
\end{Theorem}

Рассмотрим на примере.
\begin{exm}Найти управления и реакции на случайные возмущения, чтобы динамическая
система
\begin{equation*}
\begin{array}{c}
  \displaystyle\frac{dx_{1}(t)}{dt}= \Bigl(x_{1}(t)+x_{2}(t)+ e^{-t} +s_{1}(t;\mathbf{x})\Bigr)dt+ b_{1}(t;\mathbf{x})dw(t)+
  \int_{R(\gamma)}{g}_{1}(t;\mathbf{x};\gamma)\nu(dt;d\gamma), \\
   \displaystyle\frac{dx_{1}(t)}{dt}= \Bigl(x_{1}(t)x_{2}(t)+ e^{-2t}+s_{2}(t;x)\Bigr)dt+ b_{2}(t;\mathbf{x})dw(t)+
  \int_{R(\gamma)}{g}_{2}(t;\mathbf{x};\gamma)\nu(dt;d\gamma),
\end{array}
\end{equation*}
подверженная воздействию
 винеровского процесса  и совершающая скачки по действием пуассоновского процесса,
с вероятностью 1 совершала движение по поверхности
 $u(t;\mathbf{x})=x_{2}e^{-2x_{1}}$.
\end{exm}

\textbf{Решение.} Сначала построим
 систему стохастических дифференциальных уравнений,
для которой функция $u(t;\mathbf{x})=x_{2}e^{-2x_{1}}$ является
детерминированным первым интегралом.
 В соответствии с
 утверждением 2 теоремы \ref{teor-1} определим функцию, обеспечивающую
 автоморфизм функции $u(t;\mathbf{x})=x_{2}e^{-2x_{1}}$.
Тогда
\begin{equation*}
\begin{array}{cc}
  \displaystyle\frac{\partial u(\mathbf{y};t)}{\partial
  y_{1}}=-2y_{2}e^{-2y_{1}},
  & \ \ \ \
  \displaystyle\frac{\partial u(\mathbf{y};t)}{\partial
  y_{2}}=e^{-2y_{1}}
\end{array}
\end{equation*}
или
\begin{equation*}
\displaystyle\frac{\partial
\mathbf{y}(t;\mathbf{x};\gamma)}{\partial \gamma}=
\left(\begin{array}{c}
  \displaystyle\frac{\partial y_{1}(t;\mathbf{x};\gamma)}{\partial
\gamma} \\
  \displaystyle\frac{\partial y_{2}(t;\mathbf{x};\gamma)}{\partial
\gamma}\end{array}
 \right)= \left(\begin{array}{c}
         e^{-2y_{1}} \\
          2y_{2}e^{-2y_{1}}.
        \end{array}
\right)
\end{equation*}
Решение этой системы с учетом начальных данных:
\begin{equation*}
\begin{array}{cc}
  y_{1}(t; \mathbf{x};\gamma)=\displaystyle\frac{1}{2}\ln\left(2\gamma +e^{2x_{1}}\right),& \ \ \
  y_{2}(t; \mathbf{x};\gamma)=2x_{2}\gamma e^{-2x_{1}}+x_{2}.
\end{array}
\end{equation*}
Следовательно, преобразование
$g(\cdot)=(g_{1}(\cdot),g_{2}(\cdot))^{*}$, обеспечивающее функции
$u(t;\mathbf{x})=x_{2}e^{-2x_{1}}$ автоморфизм, имеет
функции-координаты:
\begin{equation*}
\begin{array}{c}
  g_{1}(t; \mathbf{x};\gamma)=\displaystyle\frac{1}{2}\ln\left(2\gamma +e^{2x_{1}}\right)-x_{1}, \ \ \ \
  g_{2}(t; \mathbf{x};\gamma)=2x_{2}\gamma e^{-2x_{1}}.
\end{array}
\end{equation*}
Теперь, в соответствии с теоремой \ref{teor-1},
строим матрицу $B$ (в данном случае -- вектор-столбец, поскольку
${\bf w}(t)$ -- одномерный винеровский процесс):
$$
B(t;{\bf x})= q_{00}\left(e^{-2x_{1}}, 2x_{2}e^{-2x_{1}}\right)^{*},
$$
где $q_{00}=q_{00}(t; \mathbf{x})$, 
$
\biggl[\dfrac{\partial B(t;{\bf x})}{\partial {\bf
x}}\biggr]=q_{00}\left(
\begin{array}{cc}
  -2e^{-2x_{1}} & 0\\
  4x_{2}e^{-2x_{1}} & 2e^{-2x_{1}}
\end{array}
\right),
$
$
\biggl[\dfrac{\partial B(t;{\bf x})}{\partial {\bf
x}}\biggr]B(t;{\bf x})=q_{00}^{2}\left(
\begin{array}{c}
 -4e^{-2x_{1}} \\
0
\end{array}
\right)=\left(
\begin{array}{c}
-4q_{00}^{2}e^{-4x_{1}} \\
0
\end{array}
\right).
$

Опираясь на формулу  (\ref{A}) строим матрицу $H(t;{\bf x})$ и вычисляем ее
определитель
$$
\det H(t;{\bf x})=\det\left(
\begin{array}{ccc}
  \vec{e}_{0} & \vec{e}_{1} & \vec{e}_{2} \\
  0 & -2x_{2}e^{-2x_{1}} & e^{-2x_{1}} \\
  f_{1} &  f_{2} &  f_{3}
\end{array}
\right)= $$
$$=
\vec{e}_{0}\Bigl(-2f_{3}x_{2} e^{-2x_{1}}-f_{2} e^{-2x_{1}}\Bigr)+
\vec{e}_{1}\Bigl(f_{1}
e^{-2x_{1}}\Bigr)+\vec{e}_{2}\Bigl(2f_{1}x_{2} e^{-2x_{1}}\Bigr),
$$
где $f_{i}=f_{i}(t; \mathbf{x})$, $i=1,2,3$. В итоге (из \eqref{A}, коэффициенты
вектора $A=A(t; \mathbf{x})$ имеют вид:
$$
\begin{array}{c}
  a_{1}= -\displaystyle\frac{f_{1}}{f_{2}+2f_{3}x_{2}}+2 q_{00}^{2}e^{-4x_{1}},\\
  a_{2}= -\displaystyle\frac{2f_{1}x_{2}}{f_{2}+2f_{3}x_{2}},
\end{array}
$$
и искомая система стохастических дифференциальных уравнений такова:
$$
\begin{array}{c}
  dx_{1}(t)=\left[-\displaystyle\frac{f_{1}}{f_{2}+2f_{3}x_{2}}+2 q_{00}^{2}e^{-4x_{1}}\right]dt + \\
   +q_{00}e^{-2x_{1}}dw(t)
  + \displaystyle\int_{R(\gamma)}\Bigl(\displaystyle\frac{1}{2}\ln\left(2\gamma +e^{2x_{1}}\right)-x_{1}\Bigr)\nu(dt;d\gamma)\\
  dx_{2}(t)=\left[-\displaystyle\frac{2f_{1}x_{2}}{f_{2}+2f_{3}x_{2}}\right]dt
  +\\
  +  q_{00}2x_{2}e^{-2x_{1}})dw(t)
  + \displaystyle\int_{R(\gamma)}(2x_{2}\gamma e^{-2x_{1}})  \nu(dt;d\gamma)
\end{array}
$$

Искомое управление -- решение системы линейных уравнений (теорема \ref{teor-2}):
$$
\begin{array}{c}
  x_{1}(t)+x_{2}(t)+ e^{-t} +s_{1}(t;\mathbf{x})=-\displaystyle\frac{f_{1}}{f_{2}+2f_{3}x_{2}}+2 q_{00}^{2}e^{-4x_{1}}, \\
  x_{1}(t)x_{2}(t)+
  e^{-2t}+s_{2}(t;\mathbf{x})=-\displaystyle\frac{2f_{1}x_{2}}{f_{2}+2f_{3}x_{2}}.
\end{array}
$$
Или
$$
\begin{array}{c}
  s_{1}(t;\mathbf{x})=-\displaystyle\frac{f_{1}}{f_{2}+2f_{3}x_{2}}+2 q_{00}^{2}e^{-4x_{1}}-x_{1}(t)-x_{2}(t)- e^{-t},  \\
  s_{2}(t;\mathbf{x})=-\displaystyle\frac{2f_{1}x_{2}}{f_{2}+2f_{3}x_{2}}-x_{1}(t)x_{2}(t)-e^{-2t},
\end{array}
$$
где $f_{i}=f_{i}(t; \mathbf{x})$, $i=1,2,3$, $q_{00}=q_{00}(t;
\mathbf{x})$. Реакция на возмущение, вызванное
винеровским процессом, определяется матрицей-столбцом с элементами
$$
\begin{array}{c}
  b_{1}(t;\mathbf{x})=q_{00}e^{-2x_{1}}, \ \ \ \
    b_{2}(t;\mathbf{x})=q_{00}2x_{2}e^{-2x_{1}}).
\end{array}
$$
Элементы для компенсатора пуассоновских скачков определяются
следующим образом:
$$\begin{array}{c}
  g_{1}(t; \mathbf{x};\gamma)=\displaystyle\frac{1}{2}\ln\left(2\gamma +e^{2x_{1}}\right)-x_{1},\\
  g_{2}(t; \mathbf{x};\gamma)=2x_{2}\gamma e^{-2x_{1}}.
\end{array}
$$

Выбор функций $f_{i}(t; \mathbf{x})$,
$i=1,2,3$ и $q_{00}(t; \mathbf{x})$ позволяет строить управление,
опираясь на какие-либо условия, например, удобством для
моделирования и реализации управления.

\begin{rem} Отметим, что можно строить программное управления
по многообразию, определяемому несколькими функциями
{\rm{\cite{D_95}}}.
\end{rem}

\section*{Выводы}
Применение теории стохастического первого интеграла позволяет
строить с вероятностью 1 программные управления для динамической
системы при наличии случайных возмущений, вызванных винеровскими {\rm{\cite{Ch-U1}}} и
пуассоновскими процессами.

\textit{Автор благодарен проф. В.А.Дубко за внимание к данной
работе.}

\end{document}